\documentstyle[12pt]{jarticle}

\begin{document}

\newtheorem{th}{THEOREM}[section]


\newtheorem{guess}{CONJECTURE}

\newtheorem{cor}[th]{COROLLARY}

\newtheorem{lm}[th]{LEMMA}

\newtheorem{rem}{Remark}

\newtheorem{ac}{Acknowledgement}

\newtheorem{ans}{ANSWER}

\renewcommand{\theequation}{\thesection.\arabic{equation}}


\title{On an Application of Hypoellipticity to Solutions of Functional
Equations}

\author{\normalsize Akira TSUTSUMI and Shigeru HARUKI}

\date{\small(June 29,1994)}

\maketitle

\begin{abstract}

 We study the regularity of solutions of functional equations of a generalized
mean value type. In this paper we give sufficient conditions for the regularity
by using hypoellipticity which is a concept of the theory of partial
differential equations. Further we also give an affirmative answer to a
conjecture of H.\'Swiatak. A part of the results was announced in the
comprehensive paper [8] on series of our joint works. To prove the regularity
of solutions of functional equation is in general one of central problem in the
theory of functional equations (see [1]).

\end{abstract}

\vspace{1cm}

\section{Introduction}

 We consider the following functional equation for the unknown $f(x)$ :

\begin{equation}
\sum_{j=1}^k
a_j(x,t)f(x+\mbox{\boldmath$\phi$}_j(t))=F(x,f(\lambda_1(x)),\cdots,f(\lambda_s(x)))+b(x,t)
\end{equation}

\vspace{0.3cm}

where \(x \in R^n\) and \(t \in R^r\). The functions in (1.1) are assumed to
satisfy the following conditions:

\vspace{0.5cm}

\begin{enumerate}

\item[A-1.]
\(a_j(x,t),  b(x,t) \in C^\infty(R^n)\) for fixed  \(t\) from an open set
$\Omega\subseteq R^r$ and \(a_j(x,t)\geq0\) in \(R^n \times \Omega,
j=1,2,\cdots,k\),

\item[A-2.]
\(a_j(x,t), b(x,t) \in C^2(R^n \times \Omega), j=1,2,\cdots,k,\)

\item[A-3.]
\(\mbox{\boldmath$\phi$}_j(t) = (\phi_{j1}(t),\cdots,\phi_{jn}(t)) \in
C^2(\Omega)\),

\item[A-4.]

\(F(x,z_1,z_2,\cdots,z_s) \in C(R^{n+s}),\)

\item[A-5.]

\(\lambda_j(x) \in C(R^n) ,j=1,2,\cdots,s.\)

\end{enumerate}

 A locally integrable function or a continuous function $f(x)$ is said to be
{\it a distribution solution} or {\it a solution in the sense of distribution}
for (1.1) if

\begin{equation}
\sum_{j=1}^k\int_{R^n}a_j(x,t)f(x+\mbox{\boldmath$\phi$}_j(t))\psi(x)dx
= \int_{R^n}
F(x,f(\lambda_1(x)),\cdots,f(\lambda_s(x)))\psi(x)dx+\int_{R^n}b(x,t)\psi(x)dx
\end{equation}

holds for every test function $\psi(x)\in {\cal D}$ and every
$t\in\Omega\subseteq R^r$, where ${\cal D}$ denotes the set of $C^{\infty}$
functions in $R^n$ with compact support \enskip (see [4],\enskip [6]p.109).\\

 The problem treated here is to investigate sufficient conditions for the
following conclusion [B] to be valid.\\

\vspace{0.8cm}

 [B] {\it Every continuous solution of (1.1) in the sense of distribution is of
class $C^\infty$, and every locally integrable solution of (1.1) in the sense
of distribution is of class $C^\infty$ almost everywhere.}

\vspace{0.8cm}

 For this problem H.\'Swiatak proved\\

\begin{th}{(\'Swiatak's Theorem([5]Theorem 3, [6]Theorem 6.1))}\\
 Assume on Eq.(1.1)\\
\\
(1) the assumptions A-1 through A-5,\\
\\
(2) There exists a $t^\circ\in\Omega$ such that
$\mbox{\boldmath$\phi$}_j(t^\circ)=0$ and $a_j(x,t^\circ)>0$ for all $x\in
R^n,j=1,2,\cdots,k,$\\
\\
(3)\quad
$\left\{\mbox{\boldmath$\phi$}^\prime_j(t^\circ)\right\},j=1,2,\cdots,k,$ spans
$R^n.$\\
\\
Then the conclusion \mbox{\rm [B]} is valid.
\end{th}

\vspace{0.5cm}

 In the statement of \'Swiatak's Theorem and in the following \quad
$\mbox{\boldmath$\phi$}_j^{(\mu)}(t)$ denotes \\
$\frac{d^\mu}{dt^\mu}\mbox{\boldmath$\phi$}_j(t),\quad
j=1,2,\cdots,k,\mu=1,2,\cdots,.$\\

 For the case of which the assumption (3) of \'Swiatak's Theorem is not
satisfied \quad \'Swiatak gave the following conjecture.

\vspace{0.5cm}

\begin{guess}{(H.\'Swiatak)}
If the assuption (3) in the Theorem 1.1 is replaced by\\
($3^{\prime}$)\quad
$\left\{\mbox{\boldmath$\phi$}{\prime\prime}_j(t^\circ)\right\},\quad
j=1,2,\cdots,k,$ \quad span $R^n,$\\
then the conclusion \mbox{\rm [B]} may be valid.
\end{guess}


Let be
$$
\mbox{\boldmath$\Psi$}(x,t)=\sum_{j=1}^k\left[\left\{2\partial_t
a_j(x,t)-\sqrt{a_j(x,t)}(\mbox{\boldmath$\phi$}_j^\prime(t)\cdot\mbox{grad}_x\sqrt{a_j(x,t)})\right\}\mbox{\boldmath$\phi$}_j^\prime(t)+a_j(x,t)\mbox{\boldmath$\phi$}_j^{\prime\prime}(t)\right]$$
where the $\cdot$ denotes the inner product of the n-dimensional real vector
space.\\
\vspace*{0.5mm}

 Our answers to the conjecture are the following.\\

\vspace*{0.5mm}

\begin{ans}{(COROLLARY 2.2)}
 If (3) in \'Swiatak's Theorem is replaced by\\
($3^{\prime\prime}$)\quad
$\left\{\mbox{\boldmath$\phi$}^\prime_j(t^\circ)\right\}$,\quad
$j=1,2,\cdots,k,$ \quad and \quad $\mbox{\boldmath$\Psi$}(x,t^\circ)$ \quad
span $R^n$ for each $x \in R^n$,\\
\vspace*{0.3cm}
then the conclusion {\rm [B]} is valid.

\end{ans}

\hspace{0.5cm}

 When we confine ourselves to the Answer 1 within \quad
$a_j(x,t),j=1,2,\cdots,k,$ \quad being positive constants, the meaning of
($3^{\prime\prime}$) in the statement of the Answer 1 is clearer.\\

\begin{ans}{(THEOREM 2.3)}
 Let $a_j(x,t)$ be identical with positive constants $a_j,
j=1,2,\cdots,k.$\quad  If the assumption (3) of \'Swiatak's Theorem is replced
by\\
($3^{\prime\prime\prime}$)\quad
$\left\{\mbox{\boldmath$\phi$}^\prime_j(t^\circ)\right\},\quad j=1,2,\cdots,k,$
\quad and \quad
$\sum_{j=1}^ka_j\mbox{\boldmath$\phi$}_j^{\prime\prime}(t^{\circ})$ span
$R^n$,\\
\smallskip
then the conclusion {\rm[B]} is valid.\\
\end{ans}

\vspace{0.5cm}

\section{Main Theorems}

 Concerning with the equation (1.1) we define the first order partial
differential operators:

\setcounter{equation}{0}
\begin{eqnarray}
L_0(x,t) & = & \mbox{\boldmath$\Psi$(x,t)} \cdot \mbox{grad}_x
\end{eqnarray}
\vspace{0.3cm}
\begin{eqnarray}
L_j(x,t) & = &
\sqrt{a_j(x,t)}\mbox{\boldmath$\phi$}_j^\prime(t)\cdot\mbox{grad}_x,\quad
j=1,\cdots,k.
\end{eqnarray}
\\
\vspace{0.5cm}
 Let $\left[L_\mu(x,t),L_\nu(x,t)\right]$ denote the commutator
$L_\mu(x,t)L_\nu(x,t)-L_\nu(x,t)L_\mu(x,t)$ of the two partial differential
operators $L_\mu(x,t)$ and $L_\nu(x,t)$ with respect to the $x$-variable.\\

\vspace*{0.5cm}
 We state the main theorems in a general form:\par
\vspace*{0.5cm}
\begin{th} \hspace{0.5cm}Assume that\\
(1) A-1 \quad through \quad A-5, \\
$(2^{\prime})$ there exists a $t^\circ\in\Omega$ such that
$\mbox{\boldmath$\phi$}_j(t^\circ)={\bf 0}$,\quad $j=1,2,\cdots,k,$\\
($3^{(4)}$) $n$ operators among the following operators are linearly
independent for each $x\in R^n$:\\
$L_{j1}(x,t^\circ),\quad \left[L_{j2}(x,t^\circ),L_{j3}(x,t^\circ)\right],\quad
\left[L_{j4}(x,t^\circ),\left[L_{j5}(x,t^\circ),L_{j6}(x,t^\circ)\right]\right],\cdots\cdots,\\
\left[L_{j\iota}(x,t^\circ),\left[L_{j\kappa}(x,t^\circ),\left[L_{j\mu}(x,t^\circ),\left[\cdots,L_{jm}(x,t^\circ)\right]\right],\cdots,\right]\right],\cdots\cdots,$\\
where each suffix \quad $j\tau$,\quad representing one integer, moves on
$\left\{0,1,2,\cdots,k\right\},$\\
span $R^n$.\\
\vspace*{0.3cm}
Then {\rm[B]} is valid.
\end{th}

\vspace{0.5cm}
\begin{cor} \hspace{0.5cm}Assume that\\
(1) A-1 \quad through \quad A-5, \\
(2) There exists a $t^\circ\in\Omega$ such that \quad
$\mbox{\boldmath$\phi$}_j(t^\circ)=${\bf 0}, \\
and \quad $a_j(x,t^\circ)>0$ \quad for all $x \in R^n$,\quad $j=1,2,\cdots,k$\\
$(3^{\prime\prime})\quad
\left\{\mbox{\boldmath$\phi$}^\prime_j(t^\circ)\right\},j=1,2,\cdots,k,$ \quad
and \quad $\mbox{\boldmath$\Psi$}(x,t^\circ)$ \quad span $R^n$ for each $x \in
R^n$,\\
\vspace*{0.3cm}
then the conclusion {\rm[B]} is valid.
\end{cor}

\vspace{0.3cm}

 Here we consider the special case of Corollary 2.2 in which\\
\begin{equation}
a_j(x,t)=a_j,\quad j=1,2,\cdots,k,
\end{equation}
where \quad $a_j$ \quad are positive constants.\\

\begin{th}\hspace{0.5cm} Assume that\\
(1)\quad A-1 \quad and \quad A-2 only for $b(x,t)$ ,and \quad A-3,A-4 \quad and
\quad A-5 ,\\
($2^{\prime\prime}$)\quad There exists a \quad $t^\circ\in\Omega$ such that
$\mbox{\boldmath$\phi$}_j(t^\circ)=0,j=1,2,\cdots,k,$ \\
($3^{\prime\prime\prime}$)\quad
$\left\{\mbox{\boldmath$\phi$}^\prime_j(t^\circ)\right\},j=1,2,\cdots,k,$ and
\quad $\sum_{j=1}^ka_j\mbox{\boldmath$\phi$}_j^{\prime\prime}(t^{\circ})$ \quad
span $R^n$,\\
\vspace*{0.3cm}
then the conclusion {\rm[B]} is valid.\\
\end{th}

\begin{rem}
 Theorem 2.1 gives a sufficient condition for {\rm [B]} even if some \quad
$a_j(x,t^\circ)$
degenerates to zero with finite order at some point $x^\circ\in R^n.$  This
shows that Theorem 2.3 is an extension of \'Swiatak's Theorem which is quoted
in the Introduction.
\end{rem}
\section{Proofs}

\setcounter{equation}{0}

\begin{lm}\hspace{0.5cm}Let $f(x)$ be a solution of $Eq.(1.1)$ in the sense of
distribution. \\Then for $t^\circ$ stated in the assumption ($2^{\prime}$) of
Theorem 2.1\hspace{0.5cm}$f(x)$ satisfies the partial differential equation:
\begin{equation}
\left[\sum_{j=1}^k(L_j^2(x,t^\circ)+\partial_t^2a_j(x,t^\circ))+L_0(x,t^\circ)\right]f(x)=\partial_t^2b(x,t^\circ)
\end{equation}
\end{lm}

\vspace{0.3cm}
{\bf Proof}.\hspace{0.2cm} Since $Eq.(1.1)$ is considered in the sense of
distribution involving a parameter $t$, we can apply $\partial_t$ on each side
of $(1.1)$ even if the solution $f(x)$ is not differentiable([6],Lemma
4.1,p.102). Thus we have
\begin{equation}
\sum_{j=1}^k\left[a_j(x,t)\mbox{\boldmath$\phi$}_j^\prime(t)\cdot\mbox{grad}_x+\partial_ta_j(x,t)\right]f(x+\mbox{\boldmath$\phi$}_j(t))=\partial_tb(x,t)
\end{equation}

and by operating $\partial_t$ further on the both side of $(3.2)$
\begin{eqnarray}
&
&\sum_{j=1}^k[a_j(x,t)(\mbox{\boldmath$\phi$}_j^\prime(t)\cdot\mbox{grad}_x)^2+2\partial_ta_j(x,t)\left(\mbox{\boldmath$\phi$}_j^\prime(t)\cdot\mbox{grad}_x\right) \nonumber \\
&
&+a_j(x,t)\left(\mbox{\boldmath$\phi$}_j^{\prime\prime}(t)\cdot\mbox{grad}_x\right)+\partial_t^2a_j(x,t)]f(x+\mbox{\boldmath$\phi$}_j(t))=\partial_t^2b(x,t)
\end{eqnarray}

is obtained.\\
Observe the following identity:
\begin{eqnarray}
\lefteqn{\left(\sqrt{a_j(x,t)}\mbox{\boldmath$\phi$}_j^\prime(t)\cdot\mbox{grad}_x\right)^2}\nonumber\\
&=  a_j(x,t)\left(\mbox{\boldmath$\phi$}_j^\prime(t)\cdot\mbox{grad}_x\right)^2
+
\sqrt{a_j(x,t)}\left(\mbox{\boldmath$\phi$}_j^\prime(t)\cdot\mbox{grad}_x\sqrt{a_j(x,t)}\right)\mbox{\boldmath$\phi$}_j^\prime(t)\cdot\mbox{grad}_x
\end{eqnarray}

Hence we have \hfill
\begin{eqnarray}
\lefteqn{a_j(x,t)\left(\mbox{\boldmath$\phi$}_j^\prime(t)\cdot\mbox{grad}_x\right)^2}\nonumber\\
&=\left(\sqrt{a_j(x,t)}\mbox{\boldmath$\phi$}_j^\prime(t)\cdot\mbox{grad}_x\right)^2
-\sqrt{a_j(x,t)}\left(\mbox{\boldmath$\phi$}_j^\prime(t)\cdot\mbox{grad}_x\sqrt{a_j(x,t)}\right)\mbox{\boldmath$\phi$}_j^\prime(t)\cdot\mbox{grad}_x.
\end{eqnarray}
Therefore each term of the principal part in $x$-derivative of $Eq.(3.3)$ is\\
\begin{equation}
\left(\sqrt{a_j(x,t)}\mbox{\boldmath$\phi$}_j^\prime(t)\cdot\mbox{grad}_x\right)^2
\end{equation}
of which sum for \quad $j=1,\cdots,k$ \quad is equal to
$\sum_{j=1}^kL_j^2(x,t)$ by $(2.2)$.\\

And the first order term of $Eq.(3.3)$ is\\
\begin{eqnarray}
\lefteqn{\sum_{j=1}^k[2\partial_ta_j(x,t)\mbox{\boldmath$\phi$}_j^\prime(t)\cdot\mbox{grad}_x} \nonumber \\
&
-\sqrt{a_j(x,t)}\left(\mbox{\boldmath$\phi$}_j^\prime(t)\cdot\mbox{grad}_x\sqrt{a_j(x,t)}\right)
\mbox{\boldmath$\phi$}_j^\prime(t)\cdot\mbox{grad}_x
+a_j(x,t)\mbox{\boldmath$\phi$}_j^{\prime\prime}(t)\cdot\mbox{grad}_x]
\end{eqnarray}

which is equal to the operator $L_0(x,t)$ defined by $(2.1)$.\\
By using these inequalities $Eq.(3.3)$ can be written in the following form:

\begin{equation}
\left[\sum_{j=1}^k(L_j^2(x,t^\circ)+\partial_t^2a_j(x,t^\circ))+L_0(x,t^\circ)\right]f(x)=\partial_t^2b(x,t^\circ)
\end{equation}

which is identical to (3.1) and thus proves Lemma.\\

\vspace{1cm}

 Let the second order partial differetial operator $P=P(x,\partial_x)$ be
written in the form:

\begin{equation}
P=\sum_{j=1}^kX_j^2+X_0+c
\end{equation}

where $X_0,X_1,\cdots,X_k$ are first order homogeneous partial differential
operators in an open set $\Lambda\subseteq R^n$ with $C^\infty$ real valued
coefficients and $c$ is a $C^\infty$ real valued function in $\Lambda$.
Namely\\
\begin{equation}
X_j=\sum_{l=1}^na_{jl}(x)\frac{\partial}{\partial x_l}
\end{equation}

where  $a_{jl}(x),\quad j=0,1,2,\cdots,k,$ \quad $l=1,2,\cdots,n$ \quad are
real valued and of class $C^\infty$.\\

\vspace{0.3cm}

 Let \quad $R[X_0,X_1,\cdots,X_k]$ \quad denote the vector space spanned by

$$X_{j1},\quad [X_{j2},X_{j3}],\quad
[X_{j4},[X_{j5},X_{j6}]],\cdots,[X_{j\iota},[X_{j\kappa},[\cdots,X_{j\mu}]]]\cdots],\cdots$$

where \quad $j\nu $\quad, representing one integer, moves over
$\left\{0,1,2,\cdots,k \right\}$.\\

 A partial differential operator $Q=Q(x,\partial_x)$ with $C^\infty$
coefficients defined in an open set $\Lambda\subseteq R^n$  is said to be {\it
hypoelliptic} in
$\Lambda$ if for any distribution solution $u(x)$ of $Qu(x)=f(x)$ in $\Lambda$
and any open subset $\Lambda^\prime$ of $\Lambda$ \quad $f(x) \in
C^\infty(\Lambda^\prime)$ implies $u(x)\in C^\infty(\Lambda^\prime)$.\\
Here the definition of a distribution solution of partial differential equation
is refered to [6] (Theorem 3.1,p.100)  and [4].
\vspace{0.5cm}

 Many sufficient conditions for hypoellipticity have been investigated. Among
them we quote the following one which is used to prove the main theorems.\\
\vspace{0.5cm}

\begin{th}{(H$\ddot{\rm o}$rmander's Theorem([2]Theorem1.11,
[3]Theorem22.2.1,p.353))}\\
Assume $R[X_0,X_1,\cdots,X_k]$ contains \quad $n$ \quad linearly independent
elements at each point $x$ $\in$ $\Lambda$.\\
Then $P=P(x,\partial_x)$ is hypoelliptic in $\Lambda$.
\end{th}
\vspace*{0.5cm}
{\bf Proof of Theorem 2.1}\quad In $Eq.(3.1)$ of Lemma 3.1 we set
$L_j(x,t^\circ)=X_j,\quad j=0,1,2,\cdots,k,$ \quad where $X_j$ are those in
$(3.10)$. By the assumption $(1)$, especially A-1 and A-2, the right hand side
$\partial_t^2b(x,t^\circ)$ is of class $C^\infty$ in $R^n$. Remark that
$c=\sum_{j=1}^k\partial_t^2a_j(x,t^\circ)$.Then the conclusion is a direct
consequence of H$\ddot{\rm o}$rmander's Theorem.\\

\vspace*{0.5cm}

{\bf Proof of Corollary 2.2} \quad The linear independence of the first order
partial differential operators $L_j(x,t^\circ)$ is equivalent to the linear
independence of $n$-vectors
$a_j(x,t^\circ)\mbox{\boldmath$\phi$}_j^\prime(t^\circ),\\
j=1,2,\cdots,k,$ \quad with the parameter $x$. Because $a_j(x,t^\circ)$ are
positive by (2) \quad the linear independence of
$a_j(x,t^\circ)\mbox{\boldmath$\phi$}_j^\prime(t^\circ),j=1,2,\cdots,k$ are
also equivalent to the linear independence of
$\mbox{\boldmath$\phi$}_j^\prime(t^\circ),j=1,2,\cdots,k$. By a similar
consideration on the first order operator $L_0(x,t^\circ)$ we can show that
assumption ($3^{\prime\prime}$) is implied in a special case which
$L_j(x,t^\circ),\quad j=0,1,2,\cdots,k,$ \quad contains n linear independent
elements at each point $x \in R^n$. This proves Corollary.\\\\
\\
{\bf Proof of Theorem 2.3} \quad Because $L_j(x,t^\circ),\quad
j=1,2,\cdots,k,\quad $ are the operators with constant coefficients \quad all
their commutators vanish. Hence the assumption $(3^{(4)})$ of Theorem 2.1 is
equivalent to the existence of $n$ linearly independent operators among the
operators $L_j(x,t^\circ),\quad j=0,1,2,\cdots,k.$ \quad  As all the
$x$-derivatives of $a_j$ are identical to zero, the assumption
($3^{\prime\prime}$) of Corollary 2.2 is equivalent to the assumption $(3''')$
of Theorem 2.3. This proves Theorem.\\

\vspace*{0.5cm}
 Thus we see that Corollary 2.2 gives directly {\it Answer 1} and Theorem 2.3
gives {\it Answer 2} to \'Swiatak's conjecture.\\
\begin{ac}
The first author wishes to express his gratitude for the support and suggestion
given by Professor Janos Acz\'el during his stay at the University of Waterloo
in 1986 and 1991. And our thanks are extended to the referee for several
remarks.\\
 This work was also partly supported by Grant-in-Aid for Scientific
Research,Grant No.62540116, and No.44146014040480, Ministry of
Education,Science and Culture.\\
\end{ac}
\end{document}